\documentclass[12pt,a4paper]{article}
\usepackage{amsmath,amssymb}
\usepackage{bm}

\makeatletter

\newcommand{\mysloppy}{\tolerance 9999 \hfuzz .5\p@ \vfuzz .5\p@}

\def\@begintheorem#1#2{\trivlist \item[\hskip \labelsep{\bf #1\ #2.}]%
	\sl
	}
\def\@opargbegintheorem#1#2#3{\trivlist
	\item[\hskip \labelsep{\bf #1\ #2\ (#3).}]%
	\sl
	}

\newcounter{myequation}[section]
\def\themyequation{(\arabic{section}.\arabic{myequation})}

\newenvironment{myequation}{%
		\global\@ignoretrue
		\refstepcounter{myequation}%
		$$}{\eqno\themyequation$$}

\makeatother

\newcommand{\ha}{Hodge algebra}
\newcommand{\cm}{Cohen-Macaulay}
\newcommand{\gor}{Gorenstein}
\newcommand{\bbm}{Buchsbaum}
\newcommand{\sr}{Stanley-Reisner}
\newcommand{\mv}{Mayer-Vietoris}

\newcommand{\card}[1]{{|#1|}}
\newcommand{\link}{{\rm link}}

\newcommand{\size}{\mathop{\rm size}}
\newcommand{\supp}{\mathop{\rm supp}}
\newcommand{\height}{\mathop{\rm ht}}
\newcommand{\depth}{\mathop{\rm depth}}

\newcommand{\Min}{\mathop{\rm Min}}

\newcommand{\Ass}{\mathop{\mbox{\rm Ass}}}

\newcommand{\core}{\mathop{\mbox{\rm core}}}
\newcommand{\geomdeg}{\mathop{\mbox{\rm geom-deg}}}
\newcommand{\arithdeg}{\mathop{\mbox{\rm arith-deg}}}

\newcommand{\ini}{\mathop{\rm in}\nolimits}
\newcommand{\ind}{\mathop{\rm ind}\nolimits}

\newcommand{\proof}{{\bf Proof.}\quad}

\newcommand{\qed}{\nolinebreak\rule{.3em}{.6em}\medskip}

\newcommand{\define}{\mathrel{:=}}

\bmdefine{\aaa}{a}
\bmdefine{\bbb}{b}
\bmdefine{\ccc}{c}
\bmdefine{\eee}{e}
\bmdefine{\ppp}{p}
\bmdefine{\xxx}{x}
\bmdefine{\yyy}{y}
\bmdefine{\zzz}{z}
\bmdefine{\NNN}{N}
\bmdefine{\ZZZ}{Z}
\bmdefine{\QQQ}{Q}
\bmdefine{\RRR}{R}
\bmdefine{\CCC}{C}

\newcommand{\mmmm}{{\frak m}}

\newcommand{\dis}{{\rm dis}}

\newcounter{thm}[section]
\def\thethm{(\arabic{section}.\arabic{thm})}

\makeatletter
\def\mynewtheorem#1#2{%
	\newenvironment{#1}{%
		\def\thethm{#2~\arabic{section}.\arabic{thm}}%
		\refstepcounter{thm}%
		\@ifnextchar[{\@@opargbegintheorem}
			{\@@begintheorem}
		}{\endtrivlist}%
}
\def\mynewtheorem#1#2{%
	\newenvironment{#1}{%
		\def\thethm{#2~\arabic{section}.\arabic{thm}}%
		\refstepcounter{thm}%
		\@ifnextchar[{\@@opargbegintheorem}
			{\@@begintheorem}
		}{\endtrivlist}%
}
\def\mynewtheoremrm#1#2{%
	\newenvironment{#1}{%
		\def\thethm{#2~\arabic{section}.\arabic{thm}}%
		\refstepcounter{thm}%
		\trivlist \item[\hskip \labelsep{\bf \thethm.}]
		}{\endtrivlist}%
}

\def\@@opargbegintheorem[#1]{%
		\trivlist \item[\hskip \labelsep{\bf \thethm\ (#1).}]\sl}
\def\@@begintheorem{%
		\trivlist \item[\hskip \labelsep{\bf \thethm.}]\sl}
\def\@@opargbegintheoremrm[#1]{%
		\trivlist \item[\hskip \labelsep{\bf \thethm\ (#1).}]\sl}
\def\@@begintheoremrm{%
		\trivlist \item[\hskip \labelsep{\bf \thethm.}]\sl}
\makeatother

\mynewtheorem{thm}{Theorem}
\mynewtheorem{cor}{Corollary}
\mynewtheorem{lemma}{Lemma}
\mynewtheorem{itemm}{Cite}
\mynewtheoremrm{example}{Example}
\mynewtheoremrm{definition}{Definition}
\mynewtheorem{claim}{Claim}
\mynewtheoremrm{remark}{Remark}
\mynewtheorem{fact}{Fact}
\mynewtheorem{prop}{Proposition}
\mynewtheorem{prob}{Problem}

\begin{document}
\mysloppy
\begin{center}
\LARGE
	{On the discrete counterparts of \cm{} algebras with straightening laws}
\end{center}
\begin{center}
\large
Mitsuhiro Miyazaki 
\end{center}
\begin{center}
Department of Mathematics,
Kyoto University of Education,
1 Fukakusa-Fujinomori-cho, Fushimi-ku, Kyoto, 612-8522 Japan
\\
E-mail:
g53448@kyokyo-u.ac.jp
\end{center}
\begin{abstract}
We study properties of a poset generating a \cm{} algebra with straightening
laws (ASL for short).
We show that if a poset $P$ generates a \cm{} ASL, then $P$ is pure
and, if $P$ is moreover \bbm{}, then $P$ is \cm.
Some results concerning a Rees algebra of an ASL defined by
a straightening closed ideal are also established.
And it is shown that if $P$ is a \cm{} poset with unique minimal element
and $Q$ is a poset ideal of $P$, then
$P\uplus Q$ is also \cm.
\end{abstract}
\medbreak\noindent
\cm{} rings,
\bbm{} rings,
Standard monomials
\bigbreak\noindent
\section{Introduction}
DeConcini, Eisenbud and Procesi
defined the notion of Hodge algebra in their article \cite{dep2}
and proved many properties of Hodge algebras.
They also showed that  many algebras appearing in algebraic geometry
and commutative ring theory
have structures of Hodge algebras.
In fact, the theory of Hodge algebras is an abstraction of 
combinatorial arguments that are used to study those rings.

A Hodge algebra is an algebra with relations which satisfy 
certain laws regulated 
by combinatorial data.
It is possible to exist many Hodge algebras with the same combinatorial
data.
And there is the simplest Hodge algebra with given combinatorial data,
called the discrete Hodge algebra.
For a given Hodge algebra, we call the discrete Hodge algebra with
the same combinatorial data the discrete counterpart of it.

Among the most important facts of 
DeConcini, Eisenbud and Procesi's 
results are
\begin{itemize}
\item
A Hodge algebra and its discrete counterpart have the same dimension.
\item
The depth of the discrete counterpart is not greater than
the depth of the original Hodge algebra.
\end{itemize}
It is known that there is a Hodge algebra whose discrete counterpart
has strictly smaller depth than the original one
\cite{hibi1}.
And we note
in Section \ref{sec bad hodge} that there is
a 
series of \cm{} Hodge algebras of dimension $n$ whose discrete
counterparts have depth 0,
where $n$ runs over the set of all 
positive integers.
So there is no hope to restrict the difference of the depth of a
Hodge algebra and that of the discrete counterpart.

But if we restrict our attention to ordinal Hodge algebras 
(algebras with straightening laws, ASL for short),
the influence of the combinatorial data to the ring theoretical
properties become greater.
So there may be a restriction to the combinatorial data by
the ring theoretical properties of an ASL.

The purpose of this article is to study properties of combinatorial data
of a \cm{} graded ASL.
Since it is equivalent to study the properties of combinatorial data
of an ASL, i.e., the properties of the partially ordered set 
(poset for short)
generating the ASL,
and to study the properties of the discrete counterpart,
our results are sometimes written in the language of posets
and sometimes in the language of commutative rings.


We have to comment the result of Terai \cite{ter}.
He asserted 
that if $A$ is a homogeneous
ASL, then the depth of its discrete counterpart is at least 
$\depth A-1$.
But his proof also works for any graded Hodge algebras
and, as is shown in Section \ref{sec bad hodge},
there is a \cm{} homogeneous Hodge algebra of dimension $n$
whose discrete counterpart is of depth 0 for any $n\geq 1$.
So we start our
 research afresh.

In Section \ref{sec bad hodge}, we give a series of examples of \cm{}
Hodge algebras of dimension $n$ whose discrete counterpart is of depth $0$,
where $n$ runs over the set of all positive integers.
It is also noted that if $n$ is an even number,
then the \ha{} of the example is \gor.

In Section \ref{sec pure}, we note that if a poset $P$ generates an
equidimensional ASL, then $P$ is pure.
In particular, if $P$ generates a \cm{} ASL, then $P$ is pure.
In Section \ref{sec bbm cm}, we show that if $P$ generates a \cm{}
ASL, and $P$ itself is \bbm, then $P$ is \cm.

Consider the following four conditions.
\begin{enumerate}
\item\label{poset int}
$P$ is a poset.
\item\label{pure poset int}
$P$ is a pure poset.
\item\label{bbm poset int}
$P$ is a \bbm{} poset.
\item\label{cm poset int}
$P$ is a \cm{} poset.
\end{enumerate}
The implications 
\ref{cm poset int}$\Rightarrow$%
\ref{bbm poset int}$\Rightarrow$%
\ref{pure poset int}$\Rightarrow$%
\ref{poset int}
are well known.
And the results of Sections \ref{sec pure} and \ref{sec bbm cm} show
that 
under the assumption that there is a \cm{} ASL generated by $P$,
\ref{bbm poset int}$\Rightarrow$%
\ref{cm poset int}
and
\ref{poset int}$\Rightarrow$%
\ref{pure poset int}
are also valid.

Section \ref{sec rees} deals with some results on Rees algebras.
If $P$ is a poset and $Q$ is a poset ideal of $P$, then a new poset
is defined by duplicationg $Q$.
See Section \ref{sec rees} for the detail.
We denote this poset by $P\uplus Q$.
If $A$ is an ASL generated by $P$,
and $Q$ satisfies a certain condition,
the straightening closed property,
then the Rees algebra of $A$ defined by $I=QA$ is an ASL
generated by $P\uplus Q$.

We show that if $P$ is \cm{} and certain conditions on reduced Euler
characteristics are satisfied, then $P\uplus Q$ is \cm.
In particular, if $P$ is a \cm{} poset with unique minimal element,
then $P\uplus Q$ is \cm{} for any poset ideal $Q$.
We also show that if $A$ is an ASL generated by $P$, $Q$ is a straightening
closed poset ideal of $P$, $P$ is \cm{} and the Rees algebra $R$
defined by $I=QA$ is \cm,
then $P\uplus Q$, 
the poset generating the ASL $R$,
is also \cm.

The author is grateful to many mathematicians,
especially the members of the seminars of Kyoto university
and Meiji university.
Sepcial thanks goes to 
Kazuhiko Kurano,
Takesi Kawasaki and
Yuji Yoshino.

\section{Preliminaries}

In this article all rings and algebras are commutative with identity.
We denote the number of elements of a finite set $X$ by $\card X$ and,
for two sets $X$ and $Y$, we denote by $X\setminus Y$ the set 
$\{x\in X|x\not\in Y\}$.
The set of integers (resp.\ non-negative integers) are denoted by
$\ZZZ$ (resp.\ $\NNN$).
Standard terminology on Hodge algebras and \sr\ rings are used freely.
See \cite{dep2}, \cite{BV}, 
\cite[Chapter II]{sta}, \cite[Chapter 5]{BH} and \cite{hoc} for example.
However, we use the term ``algebra with straightening laws''
(ASL for short) to mean an ordinal Hodge algebra.

In addition we use the following notation and convention.
\begin{itemize}
\item
We use the term poset to stand for finite partially ordered set.
\item
If $P$ is a poset, we denote the set of all the minimal 
elements of $P$ by $\min P$.
\item
If $P$ is a poset, a poset ideal of $P$ is a subset $Q$ of $P$
such that $x\in Q$, $y\in P$ and $y<x$ imply $y\in Q$.
\item
For a poset $P$, we define the order complex $\Delta(P)$ of $P$ by
\[
\Delta(P)\define\{\sigma\subseteq P\mid
\mbox{$\sigma$ is a chain}\},
\]
where a chain stands for a totally ordered subset.
We also define the reduced Euler characteristic $\tilde\chi(P)$ of
$P$ by
\[
\tilde\chi(P)\define
\tilde\chi(\Delta(P)).
\]
\item
When considering a poset,
we denote by $\infty$ (or by $-\infty$ resp.) a new
element which is larger (smaller resp.) than any other element.
\item
If $P$ is a poset and $x$, $y\in P\cup\{\infty,-\infty\}$
with $x<y$, we define
\[
(x,y)_P\define\{z\in P\mid x<z<y\}.
\]
$[x,y)_P$, $(x,y]_P$ and $[x,y]_P$ are defined similarly.
\item
We denote the \sr{} ring $k[\Delta(P)]$ by $k[P]$, 
where $k$ is a commutative ring and $P$ is a poset.
And if $k[P]$ is \cm{} (or \bbm{} resp.), 
then we say $P$ is \cm{} (\bbm{} resp.) over $k$.
\item
If $A$ is a Hodge algebra over $k$ generated by $H$ governed by $\Sigma$,
we denote by $A_{\dis}$ the discrete Hodge algebra
over $k$ generated by $H$ governed by $\Sigma$.
\item
If $X$ is a matrix with entries in a commutative ring $R$,
we denote by $I_t(X)$ the ideal of $R$ generated by all the
$t$-minors of $X$.
\item
If $B$ is an $\NNN^m$-graded ring with $B_{(0,\ldots,0)}$ a field,
then we denote by $\depth B$ the depth of $B_M$,
where $M$ is the unique $\NNN^m$-graded maximal ideal.
\end{itemize}

Next we recall the notion of a standard subset
\cite{miy3}.
\begin{definition}\label{def of sset}
Let $A$ be a Hodge algebra over $k$ generated by $H$ governed by
$\Sigma$.
A subset $\Omega$ of $H$ is called a standard subset of $H$ if for any element 
$x\in \Omega A$ and for any standard monomial $M_i$ 
appearing in the standard representation
\[
x=\sum_i b_{i}M_{i}\qquad (
\mbox{$0\neq b_{i}\in k$, $M_i$ standard})
\]
of $x$,
$\supp M_i$ meets $\Omega$.
\end{definition}
For example, a poset ideal of $H$ is a standard subset by 
\ref{dep 1.2} below.
Note that if $\Omega$ is a standard subset of $H$,
then $A/\Omega A$ is a Hodge algebra over $k$
generated by $H\setminus \Omega$ governed by $\Sigma/\Omega$.

Now we recall several facts which are used in this article.

\begin{fact}[{\cite[Theorem 6.1 and Corollary 7.2]{dep2}}]
\label{dim and depth}
If $A$ is a graded Hodge algebra over a field, then
\[
\dim A_\dis=\dim A
\]
and
\[
\depth A_\dis\leq \depth A.
\]
\end{fact}

\begin{fact}[{\cite[Proposition 1.2]{dep2}}]
\label{dep 1.2}
If $A$ is a Hodge algebra over $k$ generated by a poset
$P$ governed by $\Sigma$ and $Q$ is a poset ideal
of $P$,
then
$A/QA$ is a Hodge algebra over $k$ generated by $P\setminus Q$
governed by $\Sigma/Q$.
\end{fact}

\begin{fact}[{\cite[Proposition 5.1]{dep2}}]
\label{dep 5.1}
A square free Hodge algebra over a reduced ring is reduced.
In particular, an ASL over a field is reduced.
\end{fact}

The next result easily follows from the definition of an ASL.

\begin{fact}
\label{fact nzd}
If $A$ is an ASL generated by a poset $P$
and $x$ is the unique minimal element of $P$, then $x$ is a 
non-zero-divisor (NZD for short) of $A$.
\end{fact}

\begin{fact}[{see e.g. \cite[(5.2) Proposition]{BV}}]
\label{bv 5.2}
If $A$ is an ASL generated by a poset $P$
and $Q_1$, \ldots, $Q_n$ are poset ideals of $P$,
then
\[
(Q_1\cap \cdots\cap Q_n)A=Q_1A\cap\cdots\cap Q_nA.
\]
\end{fact}

\begin{fact}[{\cite[Proposition 1.1]{dep2}}]
\label{dep 1.1}
If $A$ is a graded Hodge algebra,
then the straightening relations give a presentation of $A$.
\end{fact}

Like \cite[II.5]{sta}, we make the following 
\begin{definition}
For a Hodge algebra $A$ over $k$ generated by $H$ governed by 
$\Sigma$, we define
\[
\vbox{\halign{\hfil$\displaystyle#$&\hfil$\displaystyle{}#{}$\hfil&
		$\displaystyle#$\hfil\cr
	\core H&\define&\bigcup_{\mbox{$N$ is a generator of $\Sigma$}}\supp N,\cr
	\core \Sigma&\define&\{\mu\in\Sigma\mid\supp\mu\subseteq\core H\},\cr
	\core A&\define&A/(H\setminus \core H)A.\cr}
	}
\]
\end{definition}
It is obvious that if $\Omega=\{x_1$, \ldots, $x_t\}$ is
a subset of $H$ such that $\Omega\cap \core H=\emptyset$,
then $\Omega$ is a standard subset of $H$
and $x_1$, \ldots, $x_t$ is an $A$-regular sequence.
In particular,
\begin{lemma}
$\core A$ is a Hodge algebra generated by $\core H$
governed by $\core \Sigma$.
Furthermore,
if $H\setminus \core H=\{x_1$, \ldots, $x_t\}$,
then $x_1$, \ldots, $x_t$ is an $A$-regular sequence
and $\core A=A/(x_1, \ldots, x_t)$.
\end{lemma}
Moreover, it is easily verified that
\[
(\core A)_\dis=\core(A_\dis).
\]
So we denote both sides by $\core A_\dis$.

\section{Examples of graded Hodge algebras whose discrete
counterparts have strictly smaller depth}
\label{sec bad hodge}

In this section, we give 
a series of examples
of Hodge algebras
whose discrete counterparts have strictly smaller depths than the 
original one.

Let $n$ be a positive integer and $X=(X_{ij})$ an $n\times n$
generic symmetric matrix, i.e.,
$X_{ij}$ with $1\leq i\leq j\leq n$ are independent indeterminates
and $X_{ji}=X_{ij}$ if $j>i$.

DeConcini and Procesi \cite[Section 5]{dp}
essentially constructed a structure of graded Hodge algebra
on
$k[X]=k[X_{ij}\mid 1\leq i\leq j\leq n]$
over $k$,
where $k$ is a commutative ring
by the following way.

Set 
\[
H\define\{
[a_1,\ldots, a_t]\mid 
\mbox{$a_i\in \NNN$ and $1\leq a_1<\cdots <a_t\leq n$}\}.
\]
For an element $\alpha=[a_1,\ldots, a_t]\in H$,
we define the size of $\alpha$ to be $t$.
And for $\alpha=[a_1,\ldots, a_t]$ and $\beta=[b_1,\ldots, b_s]\in H$,
we define the relation $\leq$ on $H$ by
\[
\alpha\leq \beta
\]
if and only if
\[
\mbox{$t\geq s$ and $a_i\leq b_i$ for $i=1$, \ldots, $s$.}
\]
It is easily verified that this relation defines a partial order on $H$.
Now set
\[
D_0\define\{[\alpha|\beta]\mid
\mbox{$\alpha$, $\beta\in H$ and $\alpha$ and $\beta$ have the same size.}\}
\]
and define the relation $<$ on $D_0$ by
\[
[\alpha|\beta]<[\gamma|\delta]
\]
if and only if 
one of the following conditions is satisfied.
\begin{itemize}
\item
$\alpha<\gamma$ in $H$.
\item
$\alpha=\gamma$ and $\beta<\delta$ in $H$.
\end{itemize}
Then it is easily verified that this relation defines a partial
order on $D_0$.

Let $\Sigma_0$ be the ideal of monomials on $D_0$
generated by
\[
\{[\alpha|\beta]\mid\alpha\not\leq\beta\}
\cup
\{[\alpha|\beta][\gamma|\delta]\mid
\mbox{$\beta\not\leq\gamma$ and $\delta\not\leq \alpha$}\}.
\]
Then
\begin{thm}[{\cite[Section 5]{dp}}]
$k[X]$ is a graded Hodge algebra over $k$ generated by
$D_0$ governed by $\Sigma_0$ with structure map
\[
D_0\ni[a_1,\ldots, a_t|b_1\ldots, b_t]\longmapsto
\det(X_{a_ib_j})\in k[X].
\]
\end{thm}

Since one can delete an element 
which is a generator of the ideal of monomials 
governing a
Hodge algebra
from the poset generating it,
if we set
\[
D\define\{[\alpha|\beta]\in D_0\mid\alpha\leq\beta\}
\]
and
\[
\Sigma\define\{\mu\in\Sigma_0\mid\supp\mu\subseteq D\},
\]
we see
\begin{cor}
$k[X]$ is a graded Hodge algebra over $k$ generated by
$D$ governed by $\Sigma$.
\end{cor}

Set
\[
\Omega_t\define\{[\alpha|\beta]\in D\mid
\size \alpha\geq t\}
\quad\mbox{and}\quad
\Omega'_t\define\{[\alpha|\beta]\in D_0\mid
\size \alpha\geq t\}.
\]
Then $\Omega_t$ 
($\Omega'_t$ resp.)
is a poset ideal of $D$ ($D_0$ resp.) and
\[
\Omega_t k[X]=
\Omega'_t k[X]=I_t(X).
\]
So by \ref{dep 1.2}, we see the following
\begin{cor}
$k[X]/I_t(X)$ is a Hodge algebra over $k$
generated by $D\setminus \Omega_t$ governed by
$\Sigma/\Omega_t$.
\end{cor}

From now on, we assume that $k$ is a field.
Then by the result of Kutz \cite{kut} (see also \cite{con}),
$k[X]/I_t(X)$ is \cm{} for any $t$.

Let us focus our attention to
the case where $t=2$ and $n\geq 3$.
Set $A\define k[X]/I_2(X)$.
Since
\[
D\setminus\Omega_2=\{[i|j]\mid
1\leq i\leq j\leq n\}
\]
and
$\Sigma/\Omega_2$ is generated by
$\{[i|j][k|l]\mid j>k$ and $l>i\}$,
the discrete counterpart $A_\dis$ of $A$ is isomorphic to
\[
k[X]/I
\]
where $I$ is the monomial ideal generated by
$
\{X_{ij}X_{kl}\mid
\mbox{$j>k$ and $l>i$}\}
$.
It is easily verified that $\sqrt I=(X_{ij}\mid 1\leq i<j\leq n)$,
so in particular
$A_\dis/\sqrt{(0)}$ is the polynomial ring in $n$ variables over $k$.
Therefore we see that
\[
\dim A=\dim A_\dis=n.
\]

In order to study the depth of $A_\dis$,
we use the technique of polarization
(see \cite[p. 107]{SV}).
Let $Y_{ij}$ ($1\leq i<j\leq n$) be $n\choose 2$ new
variables and
\[
k[X,Y]=
k[X_{ij},Y_{kl}\mid 1\leq i\leq j\leq n,
1\leq k<l\leq n]
\]
be the polynomial ring.
Then the polarization of $A_\dis$ is
\[
k[X,Y]/J,
\]
where $J$ is the monomial ideal generated by
\[
\{X_{ij}Y_{ij}\mid 1\leq i<j\leq n\}
\cup
\{X_{ij}X_{kl}\mid (i,j)\neq(k,l),\ j>k,\ l>i\}.
\]
It is known that 
$X_{12}-Y_{12}$,
$X_{13}-Y_{13}$,
\ldots,
$X_{1n}-Y_{1n}$,
$X_{23}-Y_{23}$,
\ldots,
$X_{n-1,n}-Y_{n-1,n}$
is a 
$k[X,Y]/J$-regular sequence
(see \cite[Proposition 4.3]{dep2}) and
\[
k[X,Y]/(J+(X_{ij}-Y_{ij}|1\leq i<j\leq n))
=A_\dis.
\]
So
\[
\dim A_\dis-\depth A_\dis
=
\dim(k[X,Y]/J)-\depth(k[X,Y]/J).
\]

Since $J$ is a square-free monomial ideal, 
it corresponds to a simplicial complex
by the theory of \sr{} rings.
It is easily verified that
\[
\{X_{11},X_{1n},X_{nn}\}
\cup
\{Y_{ij}\mid 1\leq i<j\leq n,\ (i,j)\neq (1,n)\}
\]
and
\[
\{X_{ij}\mid 1\leq i\leq j\leq n,\ j\leq i+1\}
\cup
\{Y_{ij}\mid 1\leq i,\ i+2\leq j\leq n\}
\]
are facets of this simplicial complex.
So it
has facets of
dimensions ${n+1\choose 2}-1$ and ${n\choose 2}+1$.
Since
\[
({n+1\choose 2}-1)-({n\choose 2}+1)=n-2,
\]
we see, by the 
theory of \sr{} rings, that
\[
\dim(k[X,Y]/J)-\depth(k[X,Y]/J)\geq n-2
\]
(see e.g. \cite[p. 370]{miy2}).
Therefore
\[
\depth A_\dis\leq 2
\]
since $\dim A_\dis=n$.

On the other hand, since no generator of $I$ involve $X_{11}$ 
and $X_{nn}$, we see that $X_{11}$, $X_{nn}$ is an 
$A_\dis$-regular sequence.
So
\[
\depth A_\dis=2.
\]
Summing up,
$A$ is a \cm{} homogeneous Hodge algebra of dimension $n$
and the discrete counterpart $A_\dis$ has depth $2$.
Moreover, $A$ is generated by $D\setminus \Omega_2$ and
$\core(D\setminus \Omega_2)=
(D\setminus \Omega_2)\setminus\{[1|1]$, $[n|n]\}$.
Therefore,
$\core A$ is a \cm{} homogeneous Hodge algebra of
dimension $n-2$ and $\depth (\core A_\dis)=0$.

\begin{remark}
It is known that $A$ is isomorphic to the second Veronese
subring of the polynomial ring in $n$ variables over $k$.
So by \cite{mok}, $A$ is \gor{} when $n$ is an
even number.
\end{remark}

\begin{remark}
As is noted above,
we can give a Hodge algebra structure on the 
polynomial ring $k[X]$.
It is also verified by the same way that the
discrete counterpart of this Hodge algebra
has strictly smaller depth than that of $k[X]$ 
if $n\geq 3$.
\end{remark}

\section{Stepping stones}
\label{sec pure}

In the following of this article, we focus our attention
to ASL and consider the following 
\begin{prob}
If there is a \cm{} ASL over $k$ generated by a poset $P$,
what can be said about $P$?
In particular, is $P$ \cm{} over $k$?
\end{prob}
To tackle this problem, we state two stepping stones and
consider the following four conditions.
\begin{enumerate}
\item\label{poset}
$P$ is a poset.
\item\label{pure poset}
$P$ is a pure poset.
\item\label{bbm poset}
$P$ is a \bbm{} poset.
\item\label{cm poset}
$P$ is a \cm{} poset.
\end{enumerate}

The implications 
\ref{cm poset}$\Rightarrow$%
\ref{bbm poset}$\Rightarrow$%
\ref{pure poset}$\Rightarrow$%
\ref{poset}
are well known.
And in the following,
we state that, 
under the assumption that there is a \cm{} ASL generated by $P$,
\ref{bbm poset}$\Rightarrow$%
\ref{cm poset}
and
\ref{poset}$\Rightarrow$%
\ref{pure poset}
are also valid.

Let $k$ be a field, $X_1$, \ldots, $X_n$ be indeterminates
over $k$ and $S=k[X_1$, \ldots, $X_n]$ a polynomial ring.
Assume that $S$ is given a graded ring structure
such that $S_0=k$ and each $X_i$ is a homogeneous element of
positive degree.

For a graded ideal $I$ of $S$,
Hartshorne \cite{har}, Sturmfels-Trung-Vogel \cite{STV}
defined the notion of geometric degree $\geomdeg I$ and
arithmetic degree $\arithdeg I$ of $I$.
By definition
\[
\deg I\leq \geomdeg I\leq \arithdeg I
\]
and
\[
\deg I=\geomdeg I\Longleftrightarrow
\mbox{$S/I$ is equidimensional,}
\]
\[
\geomdeg I=\arithdeg I\Longleftrightarrow
\mbox{$S/I$ has no embedded prime ideals.}
\]

Assume a monomial order on $S$ is settled and let $\ini(I)$
be the initial ideal of $I$ with respect to this monomial order.
It is well known that $\deg I=\deg(\ini(I))$.
And Hartshorne \cite{har}, Sturmfels-Trung-Vogel \cite{STV}
showed that
\[
\geomdeg (\ini(I))\leq\geomdeg I
\]
and
\[
\arithdeg (\ini(I))\geq \arithdeg I.
\]
So if $S/\ini(I)$ has no embedded prime ideals,
then
\[
\geomdeg(\ini(I))=\geomdeg I=\arithdeg I=\arithdeg(\ini(I))
\]
and $S/I$ is equidimensional if and only if 
$S/\ini(I)$ is equidimensional.

Assume that $A$ is a graded Hodge algebra over $k$.
Then it is well known that there is a polynomial ring $S$ with
monomial order and a graded
ideal $I$ of $S$ such that
\[
A\simeq S/I\quad\mbox{and}\quad
A_\dis\simeq S/\ini(I).
\]

Assume further that $A$ is square-free.
Then it is known
\cite[Proposition 5.1]{dep2}
that $A$ and $A_\dis$ are reduced rings.
So by the arguments above, we see the following
\begin{prop}
Let $A$ be a square-free graded \ha{} 
over a field.
Then $A$ is equidimensional if and only if 
$A_\dis$ is equidimensional.
In particular, if $A$ is a graded ASL generated by a poset $P$,
then $A$ is  equidimensional if and only if $P$ is pure.
\end{prop}
Since a \cm{} ring is equidimensional, we see the following 
\begin{cor}
Let $P$ be a poset.
If there is a \cm{} ASL generated by $P$, 
then $P$ is pure.
\end{cor}

\section{A \bbm{} poset which generate a \cm{} ASL is \cm}
\label{sec bbm cm}

In this section, we shall prove that if a poset $P$ generates a 
\cm{} ASL and if $P$ itself is \bbm{}, then
$P$ is \cm.

We begin by noting the graded version of the Theorem of 
Huckaba and Marley \cite{hm},
which is proved by noting \ref{lem 4.1} below
and reducing the proof of \cite[Proposition 3.2]{hm}
to the case of bigraded prime ideals.
\begin{thm}[The graded version of the Theorem of Huckaba-Marley]
\label{thm 4.2}
Let $A$ be a non-negatively graded Noetherian ring with 
$A_0$ a field
and $I$ a 
non-nilpotent
graded ideal of $A$.
Denote by $R$ the Rees algebra with respect to $I$ and
by $G$ the associated graded ring.
Suppose that
\[
\depth G<\depth A.
\]
Then
\[
\depth R=\depth G+1.
\]
\end{thm}
\begin{lemma}
\label{lem 4.1}
Let $R$ be an $\NNN^2$-graded ring and $I$ an 
ideal of $R$ which is homogeneous in the first grading.
If we set
\[
I^\ast\define(a\in I\mid
\mbox{$a$ is homogeneous in the second grading}),
\]
then
\[
I^\ast=(x\in I\mid
\mbox{$x$ is bihomogeneous}).
\]
In particular,
$I^\ast$ is a bigraded ideal.
\end{lemma}

Next we state the following lemma due to
Takesi Kawasaki and Kazuhiko Kurano.

\begin{lemma}
\label{lem 4.3}
Let $B$ be an $\NNN^m$-graded Noetherian ring,
$M$ a finitely generated graded $B$-module.
Then for any $i\in\ZZZ$ and $\alpha\in\ZZZ^m$,
\[
[H_N^i(M)]_\alpha
\]
is a finitely generated $B_{(0,\ldots,0)}$-module,
where $N=
\bigoplus_{\alpha\in\NNN^m\setminus\{(0,\ldots,0)\}}B_\alpha$.
\end{lemma}
\proof
We may assume that $B=B_{(0,\ldots,0)}[X_1,\ldots, X_l]$,
a polynomial ring with $\deg X_i\in\NNN^m\setminus\{(0,\ldots,0)\}$.
Then
$N=(X_1,\ldots,X_l)$
and we can compute the local cohomology with support $N$ by the 
\v Ceck complex with respect $X_1$, \ldots, $X_l$.

Since
\[
H_N^i(B)=\begin{cases}
0&\text{ if $i\neq l$}\\
X_1^{-1}\cdots X_l^{-1}B_{(0,\ldots,0)}[X_1^{-1},\ldots,X_l^{-1}]&
\text{ if $i=l$}
\end{cases}
\]
the assertion holds if $M$ is a free module.

Now we prove the assertion by the backward induction on $i$.

If $i>l$, then $H_N^i(M)=0$ since
$H_N^i(M)$ can be computed by the \v Ceck complex with 
respect to $X_1$, \ldots, $X_l$.

Now assume that $i\leq l$.
Take a short exact sequence
\begin{myequation}
\label{eq 4.1}
0\longrightarrow
K\longrightarrow
F\longrightarrow
M\longrightarrow
0
\end{myequation}
in the category of graded $B$-modules with
$F$ a free $B$-module of finite rank.
Then $[H_N^i(F)]_\alpha$ and $[H_N^{i+1}(K)]_\alpha$
are finitely generated $B_{(0,\ldots,0)}$-module
by the fact that $F$ is free and the inductive hypothesis.
So from the long exact sequence obtained by \ref{eq 4.1}, we 
see that
$[H_N^i(M)]_\alpha$ is a finitely generated 
$B_{(0,\ldots,0)}$-module.
\qed

Next we state the following
\begin{lemma}
\label{lem 4.4}
Let $k$ be an infinite field and
$G$ an $\NNN^m$-graded Hodge algebra over $k$.
Then for any $\alpha\in\ZZZ^m$,
$[H_N^i(G)]_\alpha$ is a subquotient of 
$[H_{N'}^i(G_\dis)]_\alpha$,
where $N$ (or $N'$ resp.)
is the unique $\NNN^m$-graded maximal ideal
of $G$ ($G_\dis$ resp.).
\end{lemma}
\proof
Let $P=\{x_1$, \ldots, $x_n\}$ be the poset 
which generate the Hodge algebra $G$ and $f\colon k[X_1$, \ldots, $X_n]
\rightarrow G$ be the $k$-algebra homomorphism sending $X_i$ to $x_i$,
where $k[X_1$, \ldots, $X_n]$ is the polynomial ring.
Set $\ker f=I$
and $\underline X=X_1$, \ldots, $X_n$.
Then 
\[
k[\underline X]/I\simeq G.
\]
It is well known that there is a monomial order on 
$k[\underline X]$
such that
\[
k[\underline X]/\ini(I)\simeq G_\dis.
\]

By the theory of Gr\"obner basis (see e.g. \cite[Section 15.8]{eis}),
there is an ideal $J$ of 
$k[\underline X$, $T]$ such that
\begin{itemize}
\item
$k[\underline X$, $T]/J$ is a free $k[T]$-module,
\item
$k[\underline X$, $T]/(J,T)\simeq k[\underline X]/\ini(I)$
and
\item
$k[\underline X$, $T]/(J,T-u)\simeq k[\underline X]/I$
for any $u\in k$ with $u\neq 0$,
\end{itemize}
where $T$ is a new variable.
Put 
$B=k[\underline X$, $T]$ 
and consider $B$ as an 
$\NNN^m$-graded ring by setting
$\deg T=0$ and $\deg X_i=\deg x_i$ for $i=1$, \ldots, $n$.

Then by \ref{lem 4.3}
\[
[H^i_{(\underline X)}(k[\underline X,T]/J)]_\alpha
\]
is a finitely generated $k[T]$-module for any $i\in\ZZZ$ 
and $\alpha\in\ZZZ^m$.
So we may write
\[
[H^i_{(\underline X)}(k[\underline X,T]/J)]_\alpha
=k[T]^{n_i}\oplus k[T]/(e_1)\oplus\cdots\oplus k[T]/(e_r)
\]
\[
[H^{i+1}_{(\underline X)}(k[\underline X,T]/J)]_\alpha
=k[T]^{n_{i+1}}\oplus k[T]/(f_1)\oplus\cdots\oplus k[T]/(f_s)
\]
since $k[T]$ is a principal ideal domain.

Now let $u$ be an arbitrary element of $k$.
By 
the short exact sequence
\[
0\longrightarrow
k[\underline X,T]/J\stackrel{T-u}{\longrightarrow}
k[\underline X,T]/J\longrightarrow
k[\underline X,T]/(J,T-u)\longrightarrow
0
\]
we have the following long exact sequence
\[
\vbox{\halign{&\hfil$\displaystyle#$\hfil\cr
&
[H^i_{(\underline X)}(k[\underline X,T]/J)]_\alpha
&
\stackrel{T-u}{\longrightarrow}
&
[H^i_{(\underline X)}(k[\underline X,T]/J)]_\alpha
&
\longrightarrow
[H^i_{(\underline X)}(k[\underline X,T]/(J,T-u))]_\alpha
\cr
\longrightarrow
&
[H^{i+1}_{(\underline X)}(k[\underline X,T]/J)]_\alpha
&
\stackrel{T-u}{\longrightarrow}
&
[H^{i+1}_{(\underline X)}(k[\underline X,T]/J)]_\alpha.
\cr}}
\]
So 
$[H_{(\underline X)}^i(k[\underline X,T]/(J,T-u))]_\alpha$ is
isomorphic to 
$(k[T]/(T-u))^{n_i}$,
except for finitely many $u\in k$.
And is a subquotient of the corresponding module of  any
other $u\in k$.

Since
\[
k[\underline X,T]/(J,T-u)\simeq G
\]
and so
\[
[H_N^i(G)]_\alpha\simeq
[H_{(\underline X)}^i(k[\underline X,T]/(J,T-u)]_\alpha
\]
for any $u\in k$ with $u\neq 0$,
$[H_N^i(G)]_\alpha$
is a subquotient
of 
$[H_{(\underline X)}^i(k[\underline X,T]/(J,T))]_\alpha \simeq
[H_N^i(G_\dis)]_\alpha$.
\qed

Now we state
\begin{thm}
\label{thm 4.5}
Let $A$ be a graded \cm{} 
square-free Hodge algebra over a field $k$.
Suppose that $\core A_\dis$ is \bbm{}.
Then $A_\dis$ is \cm{}.
\end{thm}
\proof
Let $H$ be the poset which generate the Hodge algebra $A$
and $\Sigma$ the ideal of monomials on $H$ which govern $A$.
Set
$
\Delta\define\{\sigma\subseteq H\mid \prod_{x_i\in\sigma}x_i\not\in\Sigma\}
$.
Then $A_\dis=k[\Delta]$.

In order to prove the theorem, we may assume,
by tensoring an infinite field containing $k$,
that $k$ is an infinite field.
And by considering $\core A$ instead of $A$,
we may assume that $A_\dis=k[\Delta]$ is \bbm{}.

We prove the theorem by induction on $\card{\ind A}$,
where $\ind A$ stands for the indiscrete part of $A$
(cf. \cite[p. 16]{dep2}).
If $\ind A=\emptyset$, then $A_\dis=A$ and the
assertion is clear.
So we assume that $\ind A\neq \emptyset$.

Take a minimal element $x$ of $\ind A$ and set $I=xA$.
Denote by $R$ the Rees algebra with respect to $I$
and by $G$ the associated graded ring.
Then $R$ is a bigraded ring and $G$ is a
bigraded Hodge algebra over $k$ 
such that $\ind G\subseteq\ind A\setminus\{x\}$
(\cite[Theorem 3.1]{dep2}).

If $G$ is \cm{}, then by the inductive hypothesis,
we see that $A_\dis=G_\dis$ is \cm{}.
So we assume that $G$ is not \cm{}.
Set $\depth G=e$ and $\dim A=d$.
And let $M$ (resp. $\mmmm$) be the unique bigraded (resp. graded)
maximal ideal of $R$ (reap. $A$).

Since $R$ and $G$ are bigraded rings,
there are two entries in the degrees of these rings.
But from now on, we use the notation and terminology concerning grading
to mean 
the grading 
newly defined by 
the Rees algebra.
Then, since $A$ is concentrated in degree 0
and $H_M^i(A)=H_\mmmm^i(A)$,
we see by the long exact sequence
of the local cohomology modules obtained by the short exact sequence
\begin{myequation}
0\longrightarrow
R_+\longrightarrow
R\longrightarrow
A\longrightarrow
0
\end{myequation}
that
\begin{myequation}
\label{eq 4.2}
[H_M^i(R_+)]_n\simeq [H_M^i(R)]_n
\end{myequation}
for any $i$, $n\in \ZZZ$ with $n\neq 0$.

On the other hand, since $IR=R_+(1)$, 
by the long exact sequence obtained by
\[
0\longrightarrow
IR\longrightarrow
R\longrightarrow
G\longrightarrow
0
\]
we see that there is an exact sequence
\begin{myequation}
\label{eq 4.3}
\vcenter{\halign{&\hfil${}#{}$\hfil\cr
	&&&\cdots&\longrightarrow&[H_M^{i-1}(G)]_n\cr
	\longrightarrow&[H_M^{i}(R_+)]_{n+1}&\longrightarrow&
		[H_M^{i}(R)]_n&\longrightarrow&[H_M^{i}(G)]_n\cr
	\longrightarrow&[H_M^{i+1}(R_+)]_{n+1}
		&\longrightarrow&\cdots\cr}}
\end{myequation}
for any $i$, $n\in\ZZZ$.

Now we recall the following result of Hochster.
\begin{thm}[{see \cite[Chapter II 4.1 Theorem]{sta}}]
\label{thm hochster}
Let $\Delta$ be a simplicial complex with vertex set 
$\{x_1$, \ldots, $x_n\}$.
Then the $\NNN^n$-graded Hilbert series  of
$H_\mmmm^i(k[\Delta])$ is
\[
\sum_{\sigma\in\Delta}
\big(\dim_k {\tilde H}^{i-\card\sigma-1}(\link_{\Delta}(\sigma);k)\big)
\prod_{x_i\in\sigma}
\frac{\lambda_i^{-1}}{1-\lambda_i^{-1}}
\]
where $\mmmm$ is the unique graded maximal ideal.
\end{thm}

We return to the proof of \ref{thm 4.5}.
By \ref{thm hochster} and \ref{lem 4.4}, we see that
\begin{myequation}
\label{eq 4.4}
[H_M^i(G)]_n=0
\end{myequation}
if $n>0$.
So we see by \ref{eq 4.3} that the map
\[
[H_M^i(R_+)]_{n+1}\longrightarrow
[H_M^i(R)]_n
\]
is an epimorphism for any $i$, $n\in\ZZZ$ with $n>0$.
On the other hand, 
\[
[H_M^i(R)]_n=0\qquad\mbox{for $n\gg 0$,}
\]
since $H_M^i(R)$ is an Artinian module.
Therefore, we see that
\[
[H_M^i(R_+)]_n\simeq[H_M^i(R)]_n=0
\]
for any $i$, $n\in\ZZZ$ with $n>0$.

So by \ref{eq 4.3}, we see that
\[
[H_M^e(G)]_0\simeq[H_M^e(R)]_0.
\]
Since $\depth R=\depth G+1=e+1$ by \ref{thm 4.2},
it follows that
\[
[H_M^e(G)]_0\simeq[H_M^e(R)]_0=0.
\]
On the other hand, $e=\depth G$ by definition,
so by \ref{eq 4.4} we see that 
there is an $n\in\ZZZ$ such that
\[
[H_M^e(G)]_n\neq0 
\quad\mbox{and}\quad
n<0.
\]
It follows from \ref{thm hochster} and \ref{lem 4.4}
that
\[
{\tilde H}^{e-2}(\link_{\Delta}(x);k)\neq0.
\]
But this 
contradicts to the assumption
that $A_\dis=k[\Delta]$ is \bbm{}.
(For characterizations of \bbm{} complexes, see
e.g. \cite{miy1}.)
\qed

\begin{remark}
Let $A$ be the Hodge algebra, of the case $n=3$, considered in 
section 
\ref{sec bad hodge}.
Then, as is easily verified,
\[
\core A_\dis
=k[X_{12}, X_{13}, X_{22}, X_{23}]/I
\]
where
\[
I=
(X_{12}^2, X_{12}X_{13}, X_{13}^2, X_{13}X_{22}, X_{13}X_{23}, X_{23}^2)
.
\]
Since
\[
I=(X_{12}, X_{13}, X_{22}, X_{23})^2\cap
(X_{12}^2, X_{13}, X_{23}^2)
\]
and
\[
I\colon
(X_{12}, X_{13}, X_{22}, X_{23})=
(X_{12}^2, X_{13}, X_{23}^2)
\]
is the primary component corresponding to the unique minimal prime
ideal $(X_{12}$, $X_{13}$, $X_{23})$ of $I$,
it is easily verified that 
\[
I\colon a=
(X_{12}^2, X_{13}, X_{23}^2)
\]
for any system of parameter $a$ of $\core A_\dis$.
So $\core A_\dis$
is \bbm{}.
And as is noted in section 
\ref{sec bad hodge},
$A$ is \cm{}.
Therefore
the square-free hypothesis in \ref{thm 4.5} is essential.
\end{remark}

\section{Rees algebras}
\label{sec rees}

As is noted below, a Rees algebra of an ASL have a structure
of an ASL under certain conditions.
In this section we study the relation between the \cm{}
property of such a Rees algebra and the property of the discrete
counterpart of it.

First we recall the definition of a straightening closed
ideal.
\begin{definition}
Let $A$ be a graded ASL over a field $k$ generated by a poset $P$
and $Q$ a poset ideal of $P$.
If every standard monomial $\mu_i$ 
appearing in the standard representation
\[
\alpha\beta=\sum_i r_i\mu_i, \qquad 0\neq r_i\in k
\]
of $\alpha\beta$ with $\alpha$, $\beta\in Q$ with $\alpha\not\sim\beta$,
has at least two factors in $Q$,
we say that $Q$ (or the ideal $QA$ of $A$) is 
straightening closed.
\end{definition}
Note that any poset ideal of a discrete ASL is 
straightening closed.

Now let $P$ be a poset and $Q$ a poset ideal of $P$.
We define the poset $P\uplus Q$ as follows
(cf. \cite[Section 9]{BV}).
Denote a copy of $Q$ by $Q^\ast$ and the element corresponding
to $x\in Q$ by $x^\ast\in Q^\ast$.
Set $P\uplus Q=P\cup Q^\ast$ as the underlying set.
And for $\alpha$, $\beta\in P\uplus Q$, we define 
$\alpha <\beta$ if and only if one of the following
three conditions is satisfied.
\begin{itemize}
\item
$\alpha$, $\beta\in P$ and $\alpha<\beta$ in $P$.
\item
$\alpha=x^\ast$, $\beta=y^\ast$ with $x$, $y\in Q$
and $x<y$ in $P$.
\item
$\alpha=x^\ast$ with $x\in Q$, $\beta\in P$
and $x\leq \beta$ in $P$.
\end{itemize}

With this notation, we recall the following fact
(see \cite[10d]{dep2} or
\cite[(9.13)]{BV}).
\begin{prop}
\label{pro 6.2}
Let $A$ be a graded ASL over a field $k$ generated by a poset $P$.
Suppose that $Q$ is a straightening closed poset ideal of $P$
and $I=QA$.
Then
\begin{enumerate}
\item
The Rees algebra $R$
with respect to $I$ is a graded ASL over $k$ generated
by $P\uplus Q$.
\item
The associated graded ring $G$ is a graded ASL
over $k$ generated by $P$
such that
$\ind G\subseteq \ind A$.
In particular, if $A$ is the discrete ASL, then so is $G$.
\end{enumerate}
\end{prop}

Now we examine the \cm{} property of $P\uplus Q$, 
where $P$ is a poset and $Q$ is a poset ideal of $P$.
\begin{thm}
\label{thm 6.3}
Let $P$ be a \cm{} poset over a field $k$ and $Q$ a poset ideal
of $P$.
If
\begin{myequation}
\tilde\chi((-\infty,x)_P)=0
\quad
\mbox{for any $x\in (P\cup\{\infty\})\setminus Q$}
\end{myequation}
then $P\uplus Q$ is also \cm{} over $k$.
\end{thm}
In order to prove the theorem,
we state several lemmas first.
\begin{lemma}
\label{lem 6.4}
A poset $P$ is \cm{} over a field $k$ if and only if
\[
{\tilde H}^i(\Delta((x,y)_P);k)=0
\]
for any $x$, $y\in P\cup\{-\infty,\infty\}$ with $x<y$ and
any $i\in\ZZZ$ with $i<\dim\Delta((x,y)_P)$.
\end{lemma}
This Lemma follows from Reisner's criterion of \cm{}
property \cite{rei} and induction.
\begin{lemma}
\label{lem 6.5}
Suppose $P$ is a poset and $Q$ is a poset ideal of $P$.
Then the natural inclusion
\[
\Delta(P)\subseteq\Delta(P\uplus Q)
\]
induces 
a 
deformation retract of the geometric
realizations.
In particular, the induced maps
\[
{\tilde H}^i(\Delta(P\uplus Q);k)\longrightarrow
{\tilde H}^i(\Delta(P);k)
\]
between the cohomologies are isomorphisms.
\end{lemma}
\proof
Let $X$ be a geometric realization of $\Delta(P)$.
Then a geometric realization of $\Delta(P\uplus Q)$
can be realized as a subset of $X\times I$,
where $I$ is the closed interval $[0,1]$,
by the following way.

Correspond $\alpha\in P\uplus Q$ with $\alpha=x^\ast$, $x\in Q$
the point $(x^\bullet,1)\in X\times I$,
where $x^\bullet$ denote the image of $x$ in $X$.
And correspond $\alpha\in P\uplus Q$ with $\alpha\in P$
the point $(\alpha^\bullet,0)\in X\times I$.
It is easily  verified that it realy 
defines a geometric realization.
Denote the geometric realization of $\Delta(P\uplus Q)$ defined
above by $Y$.
Then it is also easily verified that
$(\xi,a)\in Y$ and $0\leq b\leq a$ imply 
$(\xi,b)\in Y$, 
and so
\[
r\colon
\vtop{\halign{&\hfil$\displaystyle#$\hfil\cr
	Y&\longrightarrow &X\cr
	(\xi,a)&\longmapsto &\xi\cr}}
\]
is a retraction.
\qed

\begin{lemma}
\label{lem 6.6}
Let $P$ be a poset with unique minimal element $x_0$  and $Q$ a 
non-empty poset ideal of $P$.
Then
\[
{\tilde H}^i(\Delta((P\uplus Q)\setminus\{x_0^\ast\});k)=0
\]
for any $i\in\ZZZ$.
\end{lemma}
\proof
Set 
$\Pi=P\uplus Q\setminus\{x_0^\ast\}$,
$\Pi_1=\{\alpha\in\Pi\mid\alpha\geq x_0\}$,
$\Pi_2=\Pi\setminus\{x_0\}$
and
$\Pi_3=\Pi_1\cap \Pi_2$.
Then, since $\Pi=\Pi_1\cup\Pi_2$, we have the following \mv{}
sequence of cohomologies.
\[
\begin{array}{cccccc}
&\cdots&\longrightarrow&
{\tilde H}^{i-1}(\Delta(\Pi_1);k)\oplus
{\tilde H}^{i-1}(\Delta(\Pi_2);k)&
\longrightarrow&
{\tilde H}^{i-1}(\Delta(\Pi_3);k)\\
\longrightarrow&
{\tilde H}^{i}(\Delta(\Pi);k)&
\longrightarrow&
{\tilde H}^{i}(\Delta(\Pi_1);k)\oplus
{\tilde H}^{i}(\Delta(\Pi_2);k)&
\longrightarrow&
{\tilde H}^{i}(\Delta(\Pi_3);k)\\
\longrightarrow&\cdots
\end{array}
\]

Since
$\Pi_2=(P\setminus\{x_0\})\uplus(Q\setminus\{x_0\})$
and
$\Pi_3=P\setminus\{x_0\}$,
we see,
by \ref{lem 6.5},
that the map
\[
{\tilde H}^{j}(\Delta(\Pi_2);k)\longrightarrow
{\tilde H}^{j}(\Delta(\Pi_3);k)
\]
in the \mv{} sequence is an isomorphism for any $j$.
So we see that
\[
{\tilde H}^{i}(\Delta(\Pi);k)\simeq
{\tilde H}^{i}(\Delta(\Pi_1);k)
\]
for any $i\in\ZZZ$.

Since $\Pi_1=P$ and $P$ has a unique minimal element,
we see that
\[
{\tilde H}^{i}(\Delta(\Pi);k)\simeq
{\tilde H}^{i}(\Delta(\Pi_1);k)=0
\]
for any $i\in\ZZZ$.
\qed

Now we prove \ref{thm 6.3}.

Set $\Pi=P\uplus Q$ and we prove that 
\[
{\tilde H}^i(\Delta((\alpha,\beta)_\Pi);k)=0
\]
for any
$\alpha$, $\beta\in\Pi\cup\{-\infty,\infty\}$ with 
$\alpha<\beta$
and any
$i\in\ZZZ$ with
$i<\dim\Delta((\alpha,\beta)_\Pi)$.

The case where $\alpha\in P$ or $\beta=y^\ast$
with $y\in Q$ are clear form the  \cm{} property of $P$.

The case where $\alpha=x^\ast$ with $x\in Q$
and $\beta=y$ with $y\in P\cup\{\infty\}$.
If we put $P'=[x,y)_P$ and 
$Q'=P'\cap Q$,
then $P'$ is a poset with unique minimal element $x$ and $Q'$ is a 
poset ideal of $P'$.
Furthermore
\[
(\alpha,\beta)_\Pi=
(P'\uplus Q')\setminus\{x^\ast\}.
\]
So the result follows from \ref{lem 6.6}.

The case where $\alpha=-\infty$ and $\beta=y$ with
$y\in Q$ 
is verified by considering the poset
anti-isomorphic to $Q$ and using \ref{lem 6.6}.

The only remaining case is $\alpha=-\infty$
and $\beta=y$ with $y\in(P\cup\{\infty\})\setminus Q$.
Set $P'=(-\infty,y)_P$ and
$Q'=P'\cap Q$.
Then $Q'$ is a poset ideal of $P'$ and
\[
(\alpha,\beta)_\Pi
=P'\uplus Q'.
\]
So by \ref{lem 6.5}, we see that
\[
{\tilde H}^i(\Delta((\alpha,\beta)_\Pi);k)
\simeq
{\tilde H}^i(\Delta(P');k)
\]
for any $i\in\ZZZ$.
Since $P$ is \cm{} over $k$, we see that
\[
{\tilde H}^i(\Delta(P');k)=0
\]
for any $i\in\ZZZ$
with $i<\dim\Delta(P')$.
On the other hand,
\[
\tilde\chi(\Delta(P'))=0
\]
by assumption, so we see that
\[
{\tilde H}^i(\Delta(P');k)=0
\]
for any $i\in\ZZZ$.
Therefore 
\[
{\tilde H}^i(\Delta((\alpha,\beta)_\Pi);k)=0
\]
for any $i\in\ZZZ$.
\qed

It follows directly from \ref{thm 6.3} the following
\begin{cor}
\label{cor 6.7}
If $P$ is a \cm{} poset over a field $k$
with unique minimal element and
$Q$ is a poset ideal of $P$,
then $P\uplus Q$ is also \cm{} over $k$.
\end{cor}
Note the posets considered by
Bruns-Vetter in \cite[Section 9]{BV} are \cm{}
posets with unique minimal element.
So \ref{cor 6.7}  gives another proof
of \cite[(9.4) Theorem (b)]{BV}.

We consider the relation between the above result and \cm{}
property of Rees algebras in the sequal of this section.
First we recall the result of Trung-Ikeda
\cite{ti}.
\begin{thm}
\label{thm 6.8}
Let $(A,\mmmm)$ be a \cm{}
local ring and $I$ a proper, non-nilpotent ideal of $A$.
Denote the Rees algebra with respect to $I$ by $R$
and the associated graded ring by $G$.
Then the following are equivalent.
\begin{enumerate}
\item
$R$ is \cm{}.
\item
$G$ is \cm{} and $a(G)<0$.
\end{enumerate}
\end{thm}
(See \cite{gw} for the definition of the $a$-invariant 
$a(G)$.)
It is verified, 
by the same way as \ref{thm 4.2}, that the same conclusion
follows if one assumes that $A$ is a 
non-negatively graded ring over a field 
and $I$ is a graded ideal of $A$.
But in this case, one have to be careful enough
to interpret the $a$-invariant of $G$ to be
defined by the grading newly defined by
the
Rees algebra, not by the original grading of $A$.

Now let us examine when the $a$-invariant of $G$ is
negative, in case $A$ is a graded ASL over a field $k$
generated by a \cm{} poset $P$,
$Q$ is a straightening closed poset ideal of $P$
and $I=QA$.

By setting the original grading of $A$ as the first coordinate
and the grading of Rees algebra as the second coordinate,
we have an $\NNN^2$-grading of $G$.
It is well known that the bigraded Hilbert series
\[
H_G(\lambda ,\mu)\define
\sum_{i,j}\left(\dim_k G_{ij}\right)
\lambda^i\mu^j
\]
is a rational function, and if
\[
H_G(\lambda,\mu)=\frac{g(\lambda,\mu)}{f(\lambda,\mu)},
\]
where $f$ and $g$ are polynomials,
then,
since $G$ is \cm{},
\[
a(G)<0\Longleftrightarrow \deg_\mu g<\deg_\mu f,
\]
where $\deg_\mu$ stands for the degree with respect to $\mu$.
On the other hand,
an ASL and the discrete counterpart have the same Hilbert
function.
Therefore
\[
a(G)<0\Longleftrightarrow a(G_\dis)<0.
\]

Set
$
P=\{x_1,\ldots,x_n\}$
and
$Q=\{x_1,\ldots,x_m\}$.
Then by setting
\[
\deg x_i=
\begin{cases}
	(0,\ldots,0,1,0,\ldots,0,1)&\text{ if $i\leq m$ ($1$ in the $i$-th and the 
	last position)}\\
	(0,\ldots,0,1,0,\ldots,0,0)&\text{ if $i> m$ ($1$ in the $i$-th position)}
\end{cases}
\]
we can make $G_\dis$ an $\NNN^{n+1}$-graded ring.
The $\NNN^{n+1}$-graded Hilbert series 
$H_{G_\dis}(\lambda_1,\ldots, \lambda_n,\mu)$
of $G_\dis$ is
\[
\frac
{\sum_{\sigma\in\Delta(P)}
\left(\prod_{x_i\in\sigma\atop i\leq m}\lambda_i\mu\right)
\left(\prod_{x_i\not\in\sigma\atop i\leq m}(1-\lambda_i\mu)\right)
\left(\prod_{x_i\in\sigma\atop i>m}\lambda_i\right)
\left(\prod_{x_i\not\in\sigma\atop i>m}(1-\lambda_i)\right)}
{\prod_{i\leq m}(1-\lambda_i\mu)\prod_{i>m}(1-\lambda_i)}
\]
(see e.g. \cite[II 1.4 Theorem]{sta}).

The $\mu$-degree of the denominator is $m$ and the $\mu$-degree
of the numerator is at most $m$.
So $a(G_\dis)<0$ if and only if the coefficient of
$\mu^m$ of the numerator is zero.

The coefficient of $\mu^m$ of the numerator is
\[
\vbox{\halign{\hfil$\displaystyle{}#{}$\hfil&$\displaystyle#$\hfil\cr
	&\prod_{i=1}^m\lambda_i
		\sum_{\sigma\in\Delta(P)}(-1)^\card{Q\setminus\sigma}
	\left(\prod_{x_i\in\sigma\atop i>m}\lambda_i\right)
	\left(\prod_{x_i\not\in\sigma\atop i>m}(1-\lambda_i)\right)\cr
	=&(-1)^\card Q\prod_{i=1}^m\lambda_i
		\sum_{\tau\in\Delta(P)\atop \tau\cap Q=\emptyset}
		\left(\sum_{\upsilon\in\Delta(Q)\atop \upsilon\cup\tau\in\Delta(P)}
			(-1)^\card\upsilon\right)
		\left(\prod_{x_i\in\tau}\lambda_i\right)
		\left(\prod_{x_i\not\in\tau\atop i>m}(1-\lambda_i)\right)
			\cr
	=&(-1)^\card Q\prod_{i=1}^m\lambda_i
		\sum_{\tau\in\Delta(P)\atop \tau\cap Q=\emptyset}
		\left(\prod_{x_i\in\tau}\lambda_i\right)
		\left(\prod_{x_i\not\in\tau\atop i>m}(1-\lambda_i)\right)
		\tilde\chi(\{y\in Q\mid y<\min(\tau\cup\{\infty\})\}\cr}}
\]
Since
$
		\left(\prod_{x_i\in\tau}\lambda_i\right)
		\left(\prod_{x_i\not\in\tau,i>m}(1-\lambda_i)\right)
$
are independent polynomials,
we see that
$a(G)<0$ if and only if
$\tilde\chi(\{y\in Q\mid y<\min(\tau\cup\{\infty\})\}=0$
for any $\tau\in\Delta(P)$ with $\tau\cap Q=\emptyset$.
It is easily verified that the last condition is equivalent to
\[
\tilde\chi(\{y\in Q\mid y<x\})=0
\quad
\mbox{for any $x\in (P\cup\{\infty\})\setminus Q$}.
\]
So we see the following
\begin{lemma}
\label{lem 6.9}
In the setting above,
\[
a(G)<0
\]
if and only if
\[
\tilde\chi(\{y\in Q\mid y<x\})=0
\quad
\mbox{for any $x\in (P\cup\{\infty\})\setminus Q$}.
\]
\end{lemma}

Next we state the following
\begin{lemma}
\label{lem 6.10}
Let $P$ be a poset and $Q$ a poset ideal of $P$.
Then the following conditions are equivalent.
\begin{enumerate}
\item
\label{cond 6.10.1}
$\tilde\chi(\{y\in Q\mid y<x\})=0$
for any $x\in (P\cup\{\infty\})\setminus Q$.
\item
\label{cond 6.10.2}
$\tilde\chi((-\infty,x)_P)=0$
for any $x\in (P\cup\{\infty\})\setminus Q$.
\end{enumerate}
\end{lemma}
\proof
We first note that
\[
\vbox{%
\halign{\hfil$\displaystyle{}#{}$\hfil&$\displaystyle#$\hfil\cr
	&\tilde\chi((-\infty,x)_P)\cr
	=&\sum_{\sigma\in\Delta((-\infty,x)_P)}(-1)^{\card \sigma-1}\cr
	=&\sum_{\sigma\in\Delta((-\infty,x)_P),\sigma\subseteq Q}
		(-1)^{\card \sigma-1}
	+\sum_{\sigma\in\Delta((-\infty,x)_P),\sigma\not\subseteq Q}
		(-1)^{\card \sigma-1}\cr
	=&\tilde\chi(\{y\in Q\mid y<x\})
	+\sum_{\emptyset\neq\tau\in\Delta((-\infty,x)_P\setminus Q)}
		(-1)^{\card \tau}
		\sum_{\upsilon\in Q,\tau\cup\upsilon\in\Delta((-\infty,x)_P)}
		(-1)^{\card\upsilon-1}\cr
	=&\tilde\chi(\{y\in Q\mid y<x\})
	+\sum_{\emptyset\neq\tau\in\Delta((-\infty,x)_P\setminus Q)}
		(-1)^{\card \tau}
		\tilde\chi(\{y\in Q\mid y<\min\tau\})\cr}}
\]
for any $x\in(P\cup\{\infty\})\setminus Q$.

So \ref{cond 6.10.1}$\Longrightarrow$\ref{cond 6.10.2}
is clear and
\ref{cond 6.10.2}$\Longrightarrow$\ref{cond 6.10.1}
is proved by the
(Artinian) induction on $x$.
\qed

By \ref{lem 6.9} and \ref{lem 6.10}, we have the following
\begin{prop}
\label{prop 6.11}
Let $A$ be a graded ASL over a field $k$ generated by a \cm{}
poset $P$ and $Q$ a straightening closed poset ideal of $P$.
Denote the associated graded ring with respect to
$I=QA$ by $G$.
Then 
$a(G)<0$ if and only if 
$\tilde\chi((-\infty,x)_P)=0$ for any
$x\in(P\cup\{\infty\})\setminus Q$.
\end{prop}
Denote the Rees algebra with respect to $I$ by $R$ in the 
setting of \ref{prop 6.11}.
Then by \ref{dim and depth}, \ref{pro 6.2}, \ref{thm 6.3}
and \ref{thm 6.8},
we see the following
\begin{thm}
In the above setting, $R$ is \cm{}
if and only if 
$P\uplus Q$ is \cm{} over $k$.
In particular,
$R$ is \cm{} if and only if $R_\dis$ is \cm{}.
\end{thm}

%
%
%
%

%
%
\end{document}